# Covariance Kernels of Gaussian Markov Processes

Kerry Fendick

fendick@att.net

## Abstract

The solution to a multivariate linear Stochastic Differential Equation (SDE) with constant initial state is well known to be a Gaussian Markov process, but its covariance kernel involves the solution to an integral equation in the general case. We show that the covariance kernel has a simpler semi-parametric form for families of such solutions representing increments of a common process. We also show that a covariance kernel of a particular parametric form is necessary and sufficient for a solution to possess stationary increments and for a Gaussian process, in considerable generality, to have stationary increments and the Markov property. For a discretely sampled Gaussian process with such a parametric kernel, we derive closed-form expressions for unique maximum likelihood estimators of the parameter matrices that are unbiased, jointly sufficient, and easily computed regardless of the dimension. Using those estimators, we also derive closed-form expressions for posterior moments useful for forecasting.

# 1 Introduction

The solution to a multivariate linear Stochastic Differential Equation (SDE) with time-dependent coefficients and given initial state is a special case of a time-inhomogeneous diffusion process and an example of a Gaussian Markov process; for background, see Karatzas and Shreve [1], Sections 5.1 on pages 281-284 and Section 5.6 on pages 354- 363. The distribution of such a solution described there is well known, but its covariance kernel (also known as covariance function) is complicated in the general case, since it involves the solution to an integral equation. Its kernel is therefore not easily applied without further assumptions.

Here, we study solutions to multivariate linear SDEs representing incremental processes. Let $X_\varepsilon$ for $0 \leq \varepsilon < \delta \leq T \leq \infty$ denote a family of stochastic processes defined by

$$X_\varepsilon(t) \equiv X(t_0 + \varepsilon + t) - X(t_0 + \varepsilon) \; for \; 0 \leq t < T - \varepsilon \tag{1.1}$$

where $\{X(u): 0 \leq u < t_0 + T\}$ is some real vector-valued process of dimension $m \times 1$ and $t_0 \geq 0$ is some time value. If, for example, $X(\cdot)$ were to represent the prices of $m$ assets as a function of time, then $X_\varepsilon(\cdot)$ would represent the returns on those assets as a function of the time since $t_0 + \varepsilon$.

We will assume that each member of the family solves a linear SDE of the form

$$dX_\varepsilon(t) = \mu_\varepsilon(t)X_\varepsilon(t)dt + \sigma_\varepsilon(t)dW_\varepsilon(t) \text{ for } 0 \leq t < T - \varepsilon \tag{1.2}$$

where $\{W_\varepsilon(t): 0 \leq t < T - \varepsilon\}$ is a Brownian motion of dimension $r_\varepsilon \times 1$, $\sigma_\varepsilon(\cdot)$ is a real matrix-valued function of dimension $m \times r_\varepsilon$, and $\mu_\varepsilon(\cdot)$ is a real matrix-valued function of dimension $m \times m$. We make no assumptions about how $\mu_\varepsilon(\cdot)$, $\sigma_\varepsilon(\cdot)$, $r_\varepsilon$, and $W_\varepsilon$ depend on $\varepsilon$ otherwise. As (1.1) implies,

$$X_\varepsilon(0) = 0. \tag{1.3}$$

As Lemma 1 will show, (1.2) then implies that $EX_\varepsilon(t) = 0$ for all $0 \le t < T - \varepsilon$.

Since the positive constant $\delta$ can be arbitrarily small, one can paraphrase the above assumptions as saying that $X_0$ is defined by three properties: first, that it solves a linear SDE, second, that $X_0(t)$ represents the increment of some stochastic process over the interval $[t_0, t_0 + t]$, and third, that the first property is invariant (or stable) under perturbations to $t_0$. In the special case for which $\mu_\varepsilon(\cdot) = \mu_0(\cdot)$, $\sigma_\varepsilon(\cdot) = \sigma_0(\cdot)$, and $r_\varepsilon = r_0$ for all $0 \le \varepsilon < T$, the assumptions in (1.1) and (1.2) imply that $X_0$ is a solution to a linear SDE with stationary increments on $[0, T)$. In other words, (1.1) and (1.2) generalize the notion of stationary increments for solutions to linear SDEs. The definition of stationary increments, which is reviewed here prior to the statement of Corollary 3, is distinct from that of stationarity of the process itself.

Assuming the properties in (1.1)-(1.2) and some regulatory conditions for the coefficients, Theorem 1 and Corollary 2 here derive closed-form semi-parametric representations for the (unconditional) covariance kernel, conditional mean vector, and conditional covariance kernel of $X_0$. These representations are simpler than their well-known counterparts in Lemma 1 obtained under the weaker assumptions of (1.2) and (1.3) alone. Corollary 3 shows that a particular closed-form parametric kernel is a necessary and sufficient conditions for a solution to (1.2) and (1.3) to have stationary increments, and Corollary 5 demonstrates that univariate Gaussian processes with such a parametric kernel can exhibit the maximum possible range of correlation coefficients. Lemma 2 proves that, conditional on the state at a given time of a multivariate Gaussian process with such a parametric kernel, the distribution of the process up to that time is that of a matrix-scaled multivariate Brownian bridge ending at that state.

Although we make no explicit assumptions about the relationship between the coefficient functions $\mu_0(\cdot)$ and $\sigma_0(\cdot)$ when (1.1) and (1.2) hold, Corollary 1 shows they are in fact related by (1.1) and (1.2) themselves. That result shows that (1.1) is a non-trivial restriction on families of model satisfying (1.2) and (1.3), but also leads to the conclusion that modeling an incremental process $\{X_0(t) \equiv X(t_0 + t) - X(t_0): 0 \leq t < T\}$ as a solution to a linear SDE with coefficient functions that are not so related implicitly assumes that the suitability of the model depends on the precise choice of $t_0$.

An application for which a tractable, flexible covariance kernel is critical is Gaussian process regression, a form of machine learning; see Rasmussen and Williams [2] and Ebden, Reese, Gibson, and Aigraine [3] for background. In a common implementation of Gaussian process regression, a parametric representation for the covariance kernel is assumed, the kernel's parameters are estimated based on samples or *training data*, and the *posterior distribution* corresponding to the conditional distribution of the process with the estimated parameters is used to make predictions beyond the training data. A variety of covariance kernels have been studied, and methods for estimating their parameters have been developed typically entailing numerical algorithms and approximations.

When Gaussian process regression is applied for the prediction of multivariate time series, it is commonly called *temporal* Gaussian process regression; see Ebden, Reese, Gibson, and Aigraine [3] for a survey and Osborne, Roberts, Rogers, and Jennings [4] for further review of the literature. Appendix B of Rasmussen and Williams [2] summarizes properties of stationary Gaussian Markov processes potentially applicable for temporal Gaussian process regression. Here, we develop methods for temporal Gaussian process regression for any number of dependent time series using the parametric kernel developed under the assumption of stationary increments. Theorem 3 derive closed-form, unbiased, sufficient, maximum-likelihood estimators for the parameter matrices as a function of any discrete number of sample vectors. The estimators have the desirable properties of applying for multivariate time series with arbitrarily-spaced samples.

Temporal Gaussian process regression – like regression based on a least squares criteria – typically requires the computation of matrix inverses with a complexity of $O(n^3)$ as the number $n$ of sample vectors grows. Osborne, Roberts, Rogers, and Jennings [4] develop methods for temporal multivariate Gaussian process regression that enable the sequential updating of estimators with $O(n^2)$ complexity to include new observations or exclude old ones. Although we do not explicitly address implementation issues in the body of this paper, it is clear from the formulas for the estimators from Theorem 3 how one would sequentially update them with $O(1)$ complexity.

Theorem 4 then provides closed-form expressions for the mean and variance of the posterior distribution derived using the maximum likelihood estimators. Computing those posterior moments requires only matrix additions, matrix multiplications, and scaler divisions. Theorem 4 also shows that the mean of the posterior distribution is a linear function of time (for times beyond that of the most recent observation). Fendick [5] exploits this linearity to develop a deep machine learning algorithm for obtaining a polynomial fit of samples. In examples there for which the number of training samples is small, forecasts based on polynomials that are fit using this algorithm outperform by a wide margin forecasts based on polynomials that are fit using generalized least squares. The greater forecasting accuracy results from better use of information about correlations between the multiple variates.

Lemma 3 (used in the proof of Theorem 3) derives maximum likelihood estimators for the parameters of a discretely-sampled, matrix-scaled multivariate Brownian bridge, a result that, to the best of our knowledge, is also original. Aside from the results in this paper, we are not aware of examples of discretely sampled solutions to multivariate SDEs for which closed-form expressions have been obtained for maximum-likelihood estimators. In most cases for which the solution is not a Gaussian process, an explicit expression for the likelihood function is not known. Ait-Sahalia [6] and Li [7] have developed closed-form expansions applicable in such cases.

As the title of this paper suggests, the kernels studied here are characteristic of Gaussian Markov processes more generally. To illustrate that point for processes that need not have continuous sample paths, Corollary 4 shows, in considerable generality, that the parametric kernel described earlier is a necessary and sufficient condition for a Gaussian process to have the Markov property and stationary increments. That follows from an even more general result in Theorem 2 for wide-sense Markov processes. The theory of wide-sense Markov processes was developed in the univariate case on pages 90-91, 148-169, and 233-234 of Doob [8] and extended by Hida [9] for Gaussian Markov processes and by Beutler [10] and Mandrekar [11] for the multivariate case. The last three of those reference also characterize wide-sense Markov processes that are wide-sense stationary. Our result in Theorem 2 characterizing wide-sense Markov processes with wide-sense stationary increments is, to the best of our knowledge, new.

## 1.1 Organization of the Paper

Section 2 contains our main results characterizing the kernels of Gaussian Markov processes. It starts by reviewing solutions to multivariate, time-inhomogeneous, linear SDEs. It then describes the simplification obtained under the additional assumption of (1.1). It further characterizes those solutions for the special case of stationary increments and describes how analogous results hold for Gaussian Markov processes more generally. For the univariate case, It also describes the autocorrelation structure of solutions with stationary increments.

Section 3 shows that the parametric kernel derived in Section 2 also characterizes wide-sense Markov processes that need not be Gaussian or possess continuous sample paths.

Section 4 concludes this paper by describing how to perform Gaussian process regression using the parametric kernel derived in Sections 2 and 3 under the assumption of stationary increments. Section 4.1 derives properties of estimators for the kernel's two parameter matrices, and Section 4.2 derives expressions for the mean and variance of the posterior distribution. Since the results in Section 4 depend only on distributional properties of the solution, they depend only on the assumed zero mean and parametric kernel of the Gaussian process and not on other assumptions used in obtaining that kernel in Sections 2 and 3.

## 2    Results for Gaussian Markov Processes

Let $X_\varepsilon$ denote a random process of dimension $m \times 1$ solving (1.2) and (1.3). Implicit in (1.2) are Lipschitz and growth conditions that $\mu_\varepsilon(\cdot)$ and $\sigma_\varepsilon(\cdot)$ must satisfy; see, for example, Theorem 2.5 on page 287 of Karatzas and Shreve [1]. Later, we will replace (1.3) with the stronger assumption that (1.1) holds. But first we review well-known results obtained without that assumption.

Let

$$\Gamma_\varepsilon(t, s) \equiv E[X_\varepsilon(t)X_\varepsilon^T(s)] \; for \; 0 \leq s, t < T - \varepsilon \tag{2.1}$$

denote the covariance kernel of $X_\varepsilon$. It follows from this definition that

$$\Gamma_\varepsilon(s, t) = \Gamma_\varepsilon^T(t, s) \; for \; all \; \; 0 \leq s, t < T - \varepsilon \tag{2.2}$$

and, in particular, that

$$\Gamma_\varepsilon(t, t) = \Gamma_\varepsilon^T(t, t) \; for \; all \; 0 \leq \; t < T - \varepsilon$$

Since $\Gamma_\varepsilon(t, t)$ for each $t$ is a covariance matrix, it must be symmetric positive semi-definite. It is said to be *non-degenerate* if it is symmetric positive definite.

In the statements that follows, $I$ will denote the $m \times m$ identity matrix, and $a^{-T} \equiv (a^{-1})^T = (a^T)^{-1}$ will denote the transpose of the inverse of nonsingular square matrix $a$.

*Lemma 1: The solution to (1.2) and (1.3), where $\mu_\varepsilon(\cdot)$ and $\sigma_\varepsilon(\cdot)$ are locally bounded, is the Gaussian Markov process*

$$\left\{ X_\varepsilon(t) = f_\varepsilon(t) \int_0^t f_\varepsilon^{-1}(u)\sigma_\varepsilon(u)dW_\varepsilon(u) \colon 0 \le t < T - \varepsilon \right\}, \tag{2.3}$$

*for which*

$$\Gamma_\varepsilon(t,s) = \begin{cases} f_\varepsilon(t)h_\varepsilon(s)f_\varepsilon^T(s), & 0 \le s \le t < T - \varepsilon \\ f_\varepsilon(t)h_\varepsilon(t)f_\varepsilon^T(s), & T - \varepsilon > s > t \ge 0, \end{cases} \tag{2.4}$$

*where,*

$$f_\varepsilon(t) = I + \int_0^t \mu_\varepsilon(s)f_\varepsilon(s)ds \tag{2.5}$$

*is uniquely defined, non-singular for each t, and differentiable almost everywhere and where*

$$h_\varepsilon(s) = \int_0^s f_\varepsilon^{-1}(u)\sigma_\varepsilon(u)\sigma_\varepsilon^T(u)f_\varepsilon^{-T}(u)du. \tag{2.6}$$

*For each t, $\Gamma_\varepsilon(t,t)$ is non-degenerate if $h_\varepsilon(t)$ is non-singular.*

The solution defined by (2.3) and (2.5) was derived on page 354 of Karatzas and Shreve [1]; see also Theorem 4.10 on page 144-147 of Lipster and Shiryayev [12] and Section 5 of Beutler [10] for corresponding results obtained under variations on the assumptions here. The covariance kernel (2.4) then follows directly from (2.3), as recognized in (5.13) of Beutler [10]. The condition for $\Gamma_\varepsilon(t,t)$ to be non-degenerate was derived on pages 355-356 of Karatzas and Shreve [1].

The following is our most general result, as it applies when the coefficients of (1.2) are time inhomogeneous.

*Theorem 1: The solution to (1.1) and (1.2) for $0 \le \varepsilon < \delta \le T \le \infty$, where $\mu_0(0)$ and $\sigma_0(0)$ are nonsingular and $\mu_\varepsilon(\cdot)$ and $\sigma_\varepsilon(\cdot)$ are locally bounded, is the Gaussian Markov process*

$$\left\{ X_\varepsilon(t) = f_\varepsilon(t) \int_0^t f_\varepsilon^{-1}(u)\sigma_\varepsilon(u)dW_\varepsilon(u): \ 0 \leq t < T - \varepsilon \right\}, \tag{2.7}$$

*for which*

$$\Gamma_\varepsilon(t,s) = \begin{cases} f_\varepsilon(t)(f_\varepsilon(s)-I)\omega, & 0 \leq s \leq t < T - \varepsilon \\ (f_\varepsilon(t)-I)f_\varepsilon(s)\omega, & T - \varepsilon > s > t \geq 0, \end{cases} \tag{2.8}$$

*where $\omega$ is a nonsingular symmetric matrix of constants and where, for each t,*

$$f_\varepsilon(t) = I + \left( \int_\varepsilon^{t+\varepsilon} \sigma_0(u)\sigma_0^T(u)du \right)\omega^{-1} \tag{2.9}$$

*is non-singular. For each t, $\Gamma_\varepsilon(t,t)$ is non-degenerate if $\int_\varepsilon^{t+\varepsilon} \sigma_0(u)\sigma_0^T(u)du$ is non-singular.*

*Proof:* Under the assumptions of Theorem 1, the assumptions and conclusions of Lemma 1 hold. When (1.1) holds,

$$X_\varepsilon(t) = X_0(t+\varepsilon) - X_0(\varepsilon) \ for \ 0 \leq \varepsilon < \delta \leq T \leq \infty \ and \ 0 \leq t < T - \varepsilon. \tag{2.10}$$

Using (2.1), (2.4), and (2.10), it then follows that

$$\Gamma_\varepsilon(t,s) = \Gamma_0(t+\varepsilon, s+\varepsilon) - \Gamma_0(t+\varepsilon, \varepsilon) - \Gamma_0(\varepsilon, s+\varepsilon) + \Gamma_0(\varepsilon, \varepsilon)$$

$$= f_0(t+\varepsilon)h_0(s+\varepsilon)f_0^T(s+\varepsilon) - f_0(t+\varepsilon)h_0(\varepsilon)f_0^T(\varepsilon)$$

$$-f_0(\varepsilon)h_0(\varepsilon)f_0^T(s+\varepsilon) + f_0(\varepsilon)h_0(\varepsilon)f_0^T(\varepsilon) \ for \ s \leq t, \tag{2.11}$$

where $f_0(\cdot)$ has the properties described in Lemma 1. On the other hand, it also follows from (2.4) that

$$\Gamma_\varepsilon(t,s) = f_\varepsilon(t)f_\varepsilon^{-1}(s)\Gamma_\varepsilon(s,s) \ for \ s \leq t. \tag{2.12}$$

For the expression in (2.11) to factor into a product of the form (2.12), it is necessary that

$$\Gamma_\varepsilon(t,s) = (f_0(t+\varepsilon) - f_0(\varepsilon) + I)\left( h_0(s+\varepsilon)f_0^T(s+\varepsilon) - h_0(\varepsilon)f_0^T(\varepsilon) \right) \ for \ s \leq t \tag{2.13}$$

and that

$$f_0(\varepsilon)h_0(\varepsilon)f_0^T(s+\varepsilon) = (f_0(\varepsilon)-I)h_0(s+\varepsilon)f_0^T(s+\varepsilon) + h_0(\varepsilon)f_0^T(\varepsilon)\,. \qquad (2.14)$$

To see that (2.13) and (2.14) must hold, note that, on the right-hand side of (2.13), the product of the first summand of the first factor and the summands of the second factor accounts for all the summands in (2.11) involving $t$. By (2.5) and (2.12), the first factor of (2.13) must equal the identity matrix when $t = 0$ so that the remaining two summands of the first factor of (2.13) are required. Any additional summands in either of the two factors beyond those shown in (2.13) would either violate the property that $f_\varepsilon(0) = I$ or result in terms involving $t$ not present in (2.11). Equating the right-hand sides of (2.11) and (2.13) then yields (2.14).

Differentiating both sides of (2.14) with respect to $\varepsilon$ and setting $\varepsilon = 0$ yields

$$\dot{h}_0(0)f_0^T(s) = \dot{f}_0(0)h_0(s)f_0^T(s) + \dot{h}_0(0). \qquad (2.15)$$

By (2.5),

$$\mu_0(t) = \dot{f}_0(t)f_0^{-1}(t). \qquad (2.16)$$

Then, since $\dot{f}_0(0) = \mu_0(0)$ by (2.16) and $\dot{h}_0(0) = \sigma_0(0)\sigma_0^T(0)$ by (2.6), both are nonsingular under the assumptions; and so is $\omega \equiv \left(\dot{f}_0(0)\right)^{-1}\dot{h}_0(0)$. Rewriting (2.15) as

$$\omega - h_0(s) = \omega f_0^{-T}(s), \qquad (2.17)$$

and differentiating (2.17) with respect to $s$ then yields

$$\dot{h}_0(s) = \omega f_0^{-T}(s)\dot{f}_0^T(s)f_0^{-T}(s). \qquad (2.18)$$

Setting $s = 0$ in (2.18) and using (2.5) and the symmetry of $\dot{h}_0(0)$, we obtain

$$\omega = \dot{h}_0(0)\dot{f}_0^{-T}(0) = \omega^T, \qquad (2.19)$$

so that $\omega$ also must be symmetric. Since the left-hand side of (2.17) for given $s$ is then the difference of two symmetric matrices, the right-hand side of (2.17) also must be symmetric, so that

$$\omega f_0^{-T}(s) = f_0^{-1}(s)\omega^T = f_0^{-1}(s)\omega \, . \tag{2.20}$$

By (2.6) and (2.18),

$$\sigma_0(s)\sigma_0^T(s) = f_0(s)\omega f_0^{-T}(s)\dot{f}_0^T(s). \tag{2.21}$$

from which we conclude using (2.20) that

$$\sigma_0(s)\sigma_0^T(s) = \omega \dot{f}_0^T(s) = \dot{f}_0(s)\omega \, . \tag{2.22}$$

By (2.5) and (2.22),

$$f_0(t) = I + \left(\int_0^t \sigma_0(u)\sigma_0^T(u)du\right)\omega^{-1} \, for \, 0 \leq t < T. \tag{2.23}$$

Rearranging (2.17) and applying (2.20) again, we find that

$$h_0(s)f_0^T(s) = \omega f_0^T(s) - \omega = f_0(s)\omega - \omega \, . \tag{2.24}$$

By (2.13) and (2.24),

$$\Gamma_\varepsilon(t,s) = (f_0(t+\varepsilon) - f_0(\varepsilon) + I)\big(f_0(s+\varepsilon) - f_0(\varepsilon)\big)\omega \, \text{for} \, s \leq t \, . \tag{2.25}$$

Comparing the first factor of (2.12) and of (2.25) and using (2.23), we obtain the results in (2.8) and (2.9) for the case in which $s \leq t$. The corresponding result for the case in which $s > t$ then follows from (2.2) and (2.20).

By (2.23) and (2.24),

$$h_0(t) = \left( \int_0^t \sigma_0(u)\sigma_0^T(u)du \right) f_0^{-T}(t) .$$  (2.26)

By Lemma 1, $\Gamma_0(t,t)$ is non-degenerate if $h_0(t)$ is nonsingular, which by (2.26) and the non-singularity of $f_0(t)$ is true if $\int_0^t \sigma_0(u)\sigma_0^T(u)du$ is nonsingular. This condition for the non-degeneracy of $\Gamma_0(t,t)$ then implies the corresponding condition for the non-degeneracy of $\Gamma_\varepsilon(t,t)$ as follows from (2.8) and (2.9). ∎

Theorem 1 shows that $\Gamma_\varepsilon(\cdot,\cdot)$ has a semi-parametric representation that depends on the constant matrix $\omega$. A parametric representation would result if we were to further assume that $\sigma_0(\cdot)$ is piece-wise constant.

Our first corollary describes the relationship between $\mu_0(\cdot)$ and $\sigma_0(\cdot)$ implied by the assumptions of Theorem 1.

*Corollary 1. Under the assumptions of Theorem 1,*

$$\mu_0(t) = \sigma_0(t)\sigma_0^T(t) \left( \omega + \int_0^t \sigma_0(u)\sigma_0^T(u)du \right)^{-1} \text{ for } 0 \le t < T.$$

*Proof:* The conclusion follows immediately from (2.5) and (2.9) (or more directly from (2.16) and (2.23)).

∎

The assumptions of Theorem 1 preclude the case $\mu_0(\cdot) \equiv 0$ for which $X_0$ is a martingale. Nevertheless, if $\rho \equiv \omega^{-1}$, then

$$\mu_0(t) \to 0, \text{and } \Gamma_0(s,t) \to \int_0^{s \wedge t} \sigma_0(u)\sigma_0^T(u)du \text{ as } \rho \to 0$$

as follows from Theorem 1 and Corollary 1, consistently with the results from Lemma 1 for the martingale case.

The corollary that follows describes the conditional mean and conditional covariance kernel of $X_0$.

*Corollary 2: If, under the assumptions of Theorem 1, $\Gamma_\varepsilon(u,u)$ is non-degenerate, then*

$$E[X_\varepsilon(t)|X_\varepsilon(u)] = f_\varepsilon(t)f_\varepsilon^{-1}(u)X_\varepsilon(u) \; for \; 0 \leq u \leq t < T - \varepsilon \qquad (2.27)$$

*and*

$$E[(X_\varepsilon(t) - E[X_\varepsilon(t)|X_\varepsilon(u)])(X_\varepsilon(s) - E[X_\varepsilon(s)|X_\varepsilon(u)])^T|X_\varepsilon(u)]$$

$$= \begin{cases} f_\varepsilon(t)(f_\varepsilon^{-1}(u)f_\varepsilon(s) - I)\omega, & 0 < u \leq s \leq t < T - \varepsilon \\ (f_\varepsilon(t)f_\varepsilon^{-1}(u) - I)f_\varepsilon(s)\omega, & T - \varepsilon > s > t \geq u > 0, \end{cases} \qquad (2.28)$$

*where $\omega$ and $f_\varepsilon(\cdot)$ have the properties described by Theorem 1.*

*Proof:* By Theorem 1, $X_\varepsilon$ is Gaussian. By the results in Section 6.2.2 of Punaten and Styan [13],

$$E[X_\varepsilon(t)|X_\varepsilon(u)] = \Gamma_\varepsilon(t,u)\Gamma_\varepsilon(u,u)^- X_\varepsilon(u)$$

and

$$E[(X_\varepsilon(t) - E[X_\varepsilon(t)|X_\varepsilon(u)])(X_\varepsilon(s) - E[X_\varepsilon(s)|X_\varepsilon(u)])^T|X_\varepsilon(u)] = \Gamma_\varepsilon(t,s) - \Gamma_\varepsilon(t,u)\Gamma_\varepsilon(u,u)^-\Gamma_\varepsilon(u,s).$$

where $a^-$ denotes any generalized inverse of the matrix $a$ as defined in Section 3.5 of Petersen and Pedersen [14]. Since $\Gamma_\varepsilon(u,u)$ is assumed here to be non-degenerate, $\Gamma_\varepsilon(u,u)^- = \Gamma_\varepsilon(u,u)^{-1}$; and the statement of this corollary then follows using (2.8) and (2.20). ∎

In proving later results, we will take fuller advantage of the generality of the formulas in the above proof.

For our next corollary, $f_{i,j}^{(0)}(\cdot)$ will denote the $(i,j)^{th}$ element of $f_0(\cdot)$ for $i,j = 1, \dots, m$. We will say that $\{X_0(t): 0 \leq t < T\}$ has stationary increments if, for $s < t$, the distribution of $X_0(t) - X_0(s)$ is the same as that of $X_0(t-s)$.

*Corollary 3: If $X_0$ satisfies (1.2) and (1.3) for $0 \leq t < T$ where $\mu_0(0)$ and $\sigma_0(0)$ are nonsingular and $\mu_0(\cdot)$ and $\sigma_0(\cdot)$ are locally bounded, then $X_0$ has stationary increments if and only if*

$$\Gamma_0(t,s) = \begin{cases} s(\alpha - \beta t), & 0 \leq s \leq t < T \\ t(\alpha - \beta s), & T > s > t \geq 0, \end{cases} \qquad (2.29)$$

*where $\beta$ and $\alpha - \beta t$ for each $0 \leq t < T$ are symmetric positive semi-definite*

*Proof:* Since $X_0$ is Gaussian, it has stationary increments if and only if it has wide-sense stationary increments, i.e., if and only if

$$\Gamma_\varepsilon(t,s) = \Gamma_0(t,s) \; for \; 0 \leq \varepsilon < T \; and \; 0 \leq t < T - \varepsilon. \qquad (2.30)$$

If $X_0$ satisfying the stated assumptions has stationary increments on $[0, T)$, then (1.1) holds with $X = X_0$ and $t_0 = 0$, and the conclusions of Theorem 1 then apply. For $0 \leq t < T$, let $g(t)$ denote the matrix with $(i,j)^{th}$ element

$$g_{i,j}(t) \equiv f_{i,j}^{(0)}(t) - \left( f_{i,j}^{(0)}(0) + \dot{f}_{i,j}^{(0)}(0)t \right).$$

By the definition of a derivative,

$$g_{i,j}(t) = \text{o(t) as t} \to 0 \text{ for } i,j = 1, \dots, m. \qquad (2.31)$$

Then, using (2.9) and (2.25), we find that

$$\Gamma_\varepsilon(t,s) = \left( I + \dot{f}_0(0)t + g(t + \varepsilon) - g(\varepsilon) \right)\left( \dot{f}_0(0)s + g(s + \varepsilon) - g(\varepsilon) \right)\omega \qquad (2.32)$$

for $s \leq t$. Examining (2.32), we see that (2.30) holds if $g(\cdot) \equiv 0$. It would hold otherwise only if $g(\cdot)$ were linear, but that case is precluded by (2.31). Letting $\alpha \equiv \sigma_0(0)\sigma_0^T(0)$ and $\beta \equiv -\alpha\omega^{-1}\alpha$, the expressions in (2.29) immediately follows from (2.23) and (2.32) with $g(\cdot) \equiv 0$.

If, on the other hand, we assume that (2.29) holds, then (2.30) is satisfied, as follows from the first equality of (2.11). ∎

We now break with the assumptions of Theorem 1 to show that the parametric kernel from Corollary 3 characterizes Gaussian Markov processes with stationary increments more generally. The result follows from Theorem 4 in Section 4, as discussed following the proof of that theorem.

*Corollary 4*: Let $\{X_0(t) : 0 \leq t < T \leq \infty\}$ *denote a real Gaussian process of dimension* $m \times 1$ *with* $X_0(0) = 0$ *and with covariance kernel* $\Gamma_0(\cdot, \cdot)$ *satisfying regulatory conditions (i) and (ii) of Theorem 2. Then* $X_0$ *has stationary increments and the Markov property if and only if*

$$\Gamma_0(t, s) = \begin{cases} s(\alpha - \beta t), & 0 \leq s \leq t < T \\ t(\alpha - \beta s), & T > s > t \geq 0, \end{cases} \tag{2.33}$$

*where* $\beta$ *and* $\alpha - \beta t$ *for each* $0 \leq t < T$ *are symmetric positive definite.*

In contrast to Corollary 3, Corollary 4 does not assume that the sample path of $X_0$ is almost-surely continuous. On the other hand, Corollary 4 makes some strong assumptions about smoothness of the covariance kernel implying that $X_0$ is mean-square continuous.

We next describe properties of the covariance kernel from Corollaries 3 in the univariate case. In the statement of our results, let $R_s(t) \equiv X_0(t + s) - X_0(t)$ for $s, t \geq 0$ and $s + t < T$. For random variables $X$ and $Y$, also let $\mathrm{P}(X, Y) \equiv Cov(X, Y) / (\sqrt{Var\, X} \sqrt{Var\, Y})$ denote their correlation coefficient. As is well known, $-1 \leq \mathrm{P}(X, Y) \leq 1$ whenever it is finite.

*Proposition 1*: If $\{X_0(t) : 0 \leq t < T\}$ *with* $X_0(0) = 0$ *is a scaler-valued, continuous, zero-mean Gaussian process with the covariance kernel in (2.29), then*

$$\mathrm{P}\big(R_s(t), R_s(t + u)\big) = \frac{-s\beta}{\alpha - s\beta} \text{ and } \mathrm{P}\big(R_s^2(t), R_s^2(t + u)\big) = \frac{s^2\beta^2}{(\alpha - s\beta)^2}, \tag{2.34}$$

*for* $\geq 0$, $u \geq s \geq 0$, *and* $s + t + u < T$, *as the result of which*

$$\lim_{\beta \to -\infty} \mathrm{P}\big(R_s(t), R_s(t + u)\big) = 1, \quad \lim_{\substack{\beta \to \alpha/T \\ s \to T/2}} \mathrm{P}\big(R_s(t), R_s(t + u)\big) = -1,$$

and

$$\left(P\big(R_s(t), R_s(t+u)\big)\right)^2 \leq P\big(R_s^2(t), R_s^2(t+u)\big) \leq \left|P\big(R_s(t), R_s(t+u)\big)\right|.$$

*Proof*: Since

$$Cov\big(R_s(t), R_s(t+u)\big) = \Gamma_0(s+t, s+t+u) - \Gamma_0(t, s+t+u) - \Gamma_0(s+t, t+u) - \Gamma_0(t, t+u)$$

and

$$\sqrt{Var\, R_s(t)}\sqrt{Var\, R_s(t+u)} = \Gamma_0(s+t, s+t) - 2\,\Gamma_0(t, s+t) + \Gamma_0(t, t),$$

the first equality of (2.34) follows using (2.29). Because Gaussian moments reduce to covariances using Isselis' Theorem [15], we similarly obtain

$$Cov\big(R_s^2(t), R_s^2(t+u)\big) = 2\big(\Gamma_0(t, t+u) - \Gamma_0(t, t+s+u) - \Gamma_0(s+t, t+u) + \Gamma_0(s+t, s+t+u)\big)^2$$

and

$$\sqrt{Var\, R_s^2(t)}\sqrt{Var\, R_s^2(t+u)} = 2\big(\Gamma_0(t, t) - 2\,\Gamma_0(t, s+t) + \Gamma_0(s+t, s+t)\big)^2$$

from which the second equality of (2.34) follows also using (2.29). ∎

The quantities described in the first equality of (2.34) are commonly called the *autocorrelation coefficients of $X_0$ for time scale $s$*. The first equality of (2.34) shows that these coefficients are negative when β is positive, and positive when β is negative and that the coefficients can range over the interval $(-1, 1)$ depending on the parameters. The limit $β \to α/T$ in which the coefficients can approach $-1$ correspond to the case in which the distribution of $X_0$ approaches that of a scaled Brownian bridge on $[0, T)$.

The second equality of (2.33) implies that periods of high volatility (variability) will tend to follow periods of high volatility and periods of low volatility will tend to follow periods of low volatility. Such behavior is called volatility clustering in the financial literature.

# 3   A Generalization for Wide-Sense Markov Processes

To state our next result, we will say that $\{X_0(t): 0 \leq t \leq T < \infty\}$ is mean-square continuous if

$$E[X_0^T(t)X_0(t)] < \infty \quad \text{and} \quad \lim_{t-s\to 0} E\left[\left(X_0(t) - X_0(s)\right)^T\left(X_0(t) - X_0(s)\right)\right] = 0 \; for \; s \leq t.$$

We will also say that a process has the *wide-sense Markov property* if the best (minimum mean-square) linear predictor (BLP) of its future state given observations of its past states is always equal to the BLP given the most recent of those observations. When the BLP for a wide-sense Markov process is always equal to the most recent observation, we will call the process a *wide-sense martingale*. For formal definitions of wide-sense Markov processes and wide-sense martingales, see Definition 2.2 of Mandrekar [11].

*Theorem 2*: Let $\{X_0(t): 0 \leq t < T \leq \infty\}$ denote a real stochastic process of dimension $m \times 1$ with $X_0(0) = 0$ and assume that

  i.   *its covariance kernel $\Gamma_0(t,s)$ possesses second-order partial derivatives everywhere;*

  ii.  $\Gamma_0(t,t)$ *is an everywhere differentiable, and absolutely continuous function of $t \in [0,T)$ with a non-singular derivate at $t = 0$.*

Then, $X_0$ has wide-sense stationary increments and the wide-sense Markov property if and only if

$$\Gamma_0(t,s) = \begin{cases} s(\alpha - \beta t), & 0 \leq s \leq t < T \\ t(\alpha - \beta s), & T > s > t \geq 0, \end{cases} \tag{3.1}$$

*where $\beta$ and $\alpha - \beta t$ for each $0 \leq t < T$ are symmetric positive definite.*

*Proof*:  We first show that (3.1) is necessary. Since, by (i), $X_0$ has a continuous covariance surface, it is mean-square continuous; c.f. Section 4.1.1 of Rasmussen and Williams [2].  If we then further assume that $X_0$ is a wide-sense Markov process with $X_0(0) = 0$, Theorem 3.1 of Mandrekar [11] implies that $X_0(t) = f_0(t)U(t)$ for $0 \leq t < T$ where $f_0(t)$ is a non-singular $m \times m$ matrix for each $t$ and $f_0(0) = I$ without loss of generality and where $U$ is a wide-sense martingale of dimension $m \times 1$. Since $\Gamma_0(s,s)$ is nonsingular for each $s$ under the assumptions, so is $E[U(s)U^T(s)]$. If $U$ is any wide-sense martingale for which that last property holds, then $E[U(t)U^T(s)] = E[U(s)U^T(s)]$ for all $0 \leq s \leq t < T$ as follows from (2.12) of Beutler [10]. Therefore,

$$\Gamma_0(t,s) = E[X_0(t)X_0^T(s)] = f_0(t)E[U(t)U^T(s)]f_0^T(s) = f_0(t)h_0(s)f_0^T(s) \text{ for } 0 \leq s \leq t \qquad (3.2)$$

where $h_0(s) \equiv E[U(s)U^T(s)]$. Then, $h_0(0) = 0$ since $\Gamma_0(0,0) = 0$. And $f_0(\cdot)$ and $h_0(\cdot)$ are everywhere twice differentiable since $\Gamma_0(\cdot,s)$ and $\Gamma_0(t,\cdot)$ are assumed to be. Furthermore, $\dot{h}_0(0)$ is non-singular since $\dot{\Gamma}_0(0,0)$ is assumed to be.

Let

$$w(s,t) \equiv E\left[\left(X_0(t+s) - X_0(s)\right)^2\right]. \qquad (3.3)$$

By (3.2) and (3.3)

$$w(s,t) = f_0(t+s)h_0(t+s)f_0^T(t+s) - 2f_0(t+s)h_0(s)f_0^T(s) + f_0(s)h_0(s)f_0^T(s) \qquad (3.4)$$

for $0 \leq s \leq s+t < T$, so that $w(s,t)$ is twice differentiable in $s$. If $X_0(\cdot)$ has wide-sense stationary increments, then $\partial/\partial s \, w(s,t) = 0$ on $0 \leq s \leq s+t < T$. In particular,

$$0 = \frac{\partial}{\partial s}w(0,t) = \left(f_0(t)h_0(t)f_0^T(t)\right)' - 2f_0(t)\dot{h}_0(0) + \dot{h}_0(0) \qquad (3.5)$$

for $0 \leq t \leq T$, where the second equality of (3.5) follows from (3.2) and (3.4). Since $\dot{\Gamma}_0(t,t) =$

$= \left(f_0(t)h_0(t)f_0^T(t)\right)'$ is assumed to exist everywhere and to be integrable everywhere on $0 \leq t < T$, we conclude from (3.5) that

$$h_0(t) = f_0^{-1}(t)\left(2\left(\int_0^t f_0(u)\,du\right)\dot{h}_0(0) - \dot{h}_0(0)t\right)f_0^{-T}(t) \qquad (3.6)$$

for $0 \leq t < T$.

Likewise, wide-sense stationary increments imply that $\partial^2/\partial s^2\, w(s,t) = 0$ on $0 \leq s \leq s + t < T$. Using (3.4) and (3.6), we find that the terms of $\partial^2/\partial s^2\, w(s,t)$ involving the second derivative of $f_0(\cdot)$ vanish when $s = 0$, and we obtain

$$0 = \frac{\partial^2}{\partial s^2} w(0,t) = 2\left(\dot{f}_0(0) - \dot{f}_0(t)\right)\dot{h}_0(0)f_0^T(0)$$

for $0 \leq t < T$. Since the last two terms on the right-hand side are nonsingular, we conclude that

$$f_0(t) = f_0(0) + \dot{f}_0(0)t = I + \dot{f}_0(0)t \; for \; 0 \leq t < \text{T}. \tag{3.7}$$

From (3.2), (3.6), and (3.7), we obtain (3.1) by letting $\alpha \equiv f_0(0)\dot{h}_0(0)f_0^T(0)$ and $\beta \equiv -\dot{f}_0(0)\dot{h}_0(0)f_0^T(0)$. The stated properties of $\alpha$ and $\beta$ then follow from the properties of $\Gamma_0(t,t)$ as a non-degenerate covariance matrix.

Conversely, if (3.1) holds then $w(s,t)$ as defined in (4.7) does not vary with $s$, so that $X_0$ has wide-sense stationary increments, c.f., pages 99-101 of Doob [8]. Furthermore, $\Gamma_0(\cdot,\cdot)$ then satisfies the conditions of Theorem 2 of Beutler [10] under which $X_0$ is a wide-sense Markov process. ∎

In the case of Gaussian process, the distribution is determined by the covariance kernel, so that wide-sense stationary increments implies stationary increments, and the wide-sense Markov property implies the Markov property. Corollary 4 from Section 2 therefore follows from Theorem 2.

The functional form of the covariance kernel obtained in (4.5) for wide-sense Markov processes is the same as that in (2.4) for solutions to linear SDEs. An expression of the form in (2.6) therefore can be justified for a wide-sense Markov processes with a sufficiently well-behaved covariance function. Consequently, Theorem 1 can be similarly recast with minor modification to its proof to describe Gaussian Markov processes with well-behaved covariance functions. In that case, the function $\sigma_0(\cdot)\sigma_0^T(\cdot)$ appearing in (2.9) is naturally described as a time-varying instantaneous covariance matrix.

# 4 Gaussian Markov Regression

The relative simplicity of the parametric covariance kernel derived for the case of stationary increments leads to simple formulas for parameter estimation and prediction of multivariate time series.

## 4.1 Parameter Estimation

Theorem 3 is our main result on parameter estimation.

*Theorem 3: If $\{X(t): 0 \leq t < T \leq \infty\}$ is a zero-mean Gaussian process of dimension $m \times 1$ with $X(0) = 0$ and the covariance kernel*

$$\Gamma(t,s) = \begin{cases} s(\alpha - \beta t), & 0 \leq s \leq t < T \\ t(\alpha - \beta s), & T > s > t \geq 0, \end{cases} \tag{4.1}$$

*where $\beta$ and $\alpha - \beta t$ for each $0 < t \leq T$ are symmetric positive definite for each $0 < t < T$, and if samples $X(t_i) = x_i$ are known for $0 < t_1 < t_2 < \ldots < t_n < T$, then*

$$\hat{\alpha} = \frac{1}{n-1} \sum_{i=1}^{n-1} \frac{(t_{i+1}x_i - t_i x_{i+1})(t_{i+1}x_i - t_i x_{i+1})^T}{t_i t_{i+1}(t_{i+1} - t_i)} \tag{4.2}$$

*and*

$$\hat{\beta} = \frac{\hat{\alpha}}{t_n} - \frac{x_n x_n^T}{t_n^2} \tag{4.3}$$

*are the unique maximum-likelihood, jointly sufficient, unbiased estimators of the parameter matrices.*

Our proof will depend on the following two lemmas. The first describes a useful property of Gaussian processes with the kernel in (4.1).

*Lemma 2: If $\{X(t): 0 \leq t < T \leq \infty\}$ is a zero-mean Gaussian process of dimension $m \times 1$ with covariance kernel from (4.1), and if $X(u)$ is given for some $0 < u < T$, then $\{X(t): 0 \leq t \leq u\}$ is conditionally a Gaussian process with*

$$E[X(t)|X(u)] = tu^{-1}X(u) \text{ for } 0 \le t \le u$$

*and*

$$E[(X(t) - E[X(t)|X(u)])(X(s) - E[X(s)|X(u)])^T|X(u)]$$

$$= \begin{cases} (I - tu^{-1}I)s\alpha, & 0 \le s \le t \le u, \\ t\alpha(I - su^{-1}I), & u \ge s > t \ge 0. \end{cases}$$

*Proof:* By definition, a process is Gaussian if its finite-dimensional distributions are Gaussian. For any set of random variables with a joint Gaussian distribution, it is well known that the conditional distribution of any subset - given the values of the others – is Gaussian, c.f., page 522 of Rao [16]. Therefore $\{X(t): 0 \le t \le u < T\}$ is conditionally a Gaussian process given $X(u)$.

As in the proof of Corollary 2,

$$E[X(t)|X(u)] = \Gamma(t,u)\Gamma(u,u)^- X(u) = tu^{-1}\Gamma(u,u)\Gamma(u,u)^- X(u) = tu^{-1}X(u)$$

for $0 \le t \le u < T$ where the second equality follows from (4.1) and the last by Definition 2.1 of a generalized inverse in Rao and Mitra [17]. (The product $\Gamma(u,u)\Gamma(u,u)^-$ is not in general equal to the identity matrix, but it has the same effect when pre-multiplying any vector.) As in the proof of Corollary 2,

$$E[(X(t) - E[X(t)|X(u)])(X(s) - E[X(s)|X(u)])^T|X(u)] = \Gamma(t,s) - \Gamma(t,u)\Gamma_\varepsilon(u,u)^-\Gamma(u,s)$$

$$= \Gamma(t,s) - stu^{-2}\Gamma(u,u)\Gamma(u,u)^-\Gamma(u,u) = \Gamma(t,s) - stu^{-2}\Gamma(u,u) = (I - tu^{-1}I)s\alpha$$

where the second equality follows from (4.1), the third from Definition 2.2 of a generalized inverse in Rao and Mitra [17] (shown there to be equivalent to their Definition 2.1 used above), and the last by (3.1) again. ∎

We recognize the conditional moments derived in Lemma 2 as those of a matrix-scaled, multivariate Brownian bridge ending at the conditioned state. The next lemma describes the maximum likelihood estimator of its scaling matrix with discrete sampling.

Lemma *3: If $\{B(t): 0 \leq t \leq t_n < \infty\}$ is a Gaussian process of dimension $m \times 1$ with*

$$E[B(t)] = t t_n^{-1} x_n \ \ for \ 0 \leq t \leq t_n$$

*and*

$$E[B(t)B^T(t)] = \begin{cases} (I - t t_n^{-1} I) s \alpha, & 0 \leq s \leq t \leq t_n, \\ t \alpha (I - s t_n^{-1} I), & t_n > s > t \geq 0. \end{cases}$$

*where $\alpha$ is symmetric positive definite and if samples $B(t_i) = x_i$ are known for $0 < t_1 < t_2 < \ldots < t_{n-1} < T$, then (4.2) is both a sufficient estimator and the unique maximum likelihood estimator of $\alpha$.*

*Proof:* The random vector of dimension $m(n-1) \times 1$ that is constructed by concatenating $B(t_1), B(t_2), \ldots, B(t_{n-1})$ has a Gaussian distribution. By Lemma 2, the assumptions imply that its mean is equal to $vec(M)$, where M is the matrix with $j^{th}$ column equal to $t_j t_n^{-1} x_n$ for $j = 1, \ldots, n-1$, and its covariance matrix is equal to the Kronecker product $T \otimes \alpha$, where $T$ is the symmetric matrix of dimension $(n-1) \times (n-1)$ with $(i,j)th$ element equal to $t_i(1 - t_j / t_n)$ for $i \leq j$; see Section 10.2 of Petersen and Pedersen [14] for definitions of the vec operator and Kronecker product. Let $Y \equiv X - M$, where $X$ is the matrix with $j^{th}$ column equal to $x_j$ for $j = 1, \ldots, n-1$, and $y \equiv vec(Y)$. If

$$p_1(x_1, x_2, \ldots, x_{n-1} | x_n) \equiv \frac{P(B(t_1) \in dx_1, \ldots, B(t_{n-1}) \in dx_{n-1})}{dx_1 \ldots dx_{n-1}} \tag{4.4}$$

then its log-likelihood function is

$$\mathcal{L}_1 = const - \frac{1}{2}\log(|T \otimes \alpha|) - \frac{1}{2} y^T (T^{-1} \otimes \alpha^{-1}) y \tag{4.5}$$

where $|\cdot|$ denotes the matrix determinate operator. (The notation on the left-hand side of (4.4) is suggestive of our later use of this lemma.) Let

$$\alpha = \sum_{j=1}^{m} v_j v_j^T \tag{4.6}$$

denote a rank-1 decomposition of the matrix $\alpha$ in which each of the $v_i's$ is a vector of dimension $m \times 1$.

Also let $z \equiv (T^{-1} \otimes \alpha^{-1})y = (T^{-1} \otimes \alpha^{-1})vec(Y)$. By (520) of Peterson and Pedersen [14],

$$z = vec(\alpha^{-1}YT^{-1}). \qquad (4.7)$$

Then,

$$y^T(T^{-1} \otimes \alpha^{-1})y = z^T(T \otimes \alpha)z = \sum_{j=1}^{m} z^T\left(T \otimes \left(v_j v_j^T\right)\right)z$$

$$= \sum_{j=1}^{m} v_j^T \alpha^{-1}YT^{-1}TT^{-1}Y^T\alpha^{-1}v_j$$

$$= tr\left(\sum_{i=j}^{m} \alpha^{-1}YT^{-1}Y^T\alpha^{-1}v_j v_j^T\right)$$

$$= tr(\alpha^{-1}YT^{-1}Y^T), \qquad (4.8)$$

where the first equality follows from the definition of $z$, the second from (4.6) here and (506) of Petersen and Pederson [14], the third from (4.7) here and (524) of Petersen and Pedersen [14], the forth from well-known properties of the matrix trace, and the fifth from (4.6) again.

By (4.5) and (4.8),

$$\mathcal{L}_1 = const - \frac{1}{2}\log(|T \otimes \alpha|) - \frac{1}{2}tr(\alpha^{-1}YT^{-1}Y^T). \qquad (4.9)$$

from which we obtain, for each $1 \le i \le j \le m$,

$$\frac{\partial}{\partial \alpha_{i,j}}\mathcal{L}_1 = -\frac{1}{2}tr\left((T^{-1} \otimes \alpha^{-1})\left(T \otimes \frac{\partial \alpha}{\partial \alpha_{i,j}}\right)\right) + \frac{1}{2}tr\left(\alpha^{-1}YT^{-1}Y^T\alpha^{-1}\frac{\partial \alpha}{\partial \alpha_{i,j}}\right)$$

$$= -\frac{1}{2}tr\left(I_{(n-1)\times(n-1)} \otimes \alpha^{-1}\frac{\partial \alpha}{\partial \alpha_{i,j}}\right) + \frac{1}{2}tr\left(\alpha^{-1}YT^{-1}Y^T\alpha^{-1}\frac{\partial \alpha}{\partial \alpha_{i,j}}\right)$$

$$= \frac{1}{2}tr\left(-(n-1)\alpha^{-1}\frac{\partial \alpha}{\partial \alpha_{i,j}}\right) + \frac{1}{2}tr\left(\alpha^{-1}YT^{-1}Y^T\alpha^{-1}\frac{\partial \alpha}{\partial \alpha_{i,j}}\right) \qquad (4.10)$$

for $1 \leq i \leq j \leq m$ where $I_{(n-1)\times(n-1)}$ denotes the identity matrix of dimension $(n-1) \times (n-1)$. The first equality of (4.10) follows by applying (57), (512), and (137) of Petersen and Pederson [14] to the first non-constant summand on the right-hand side of (4.9) and by applying (124) and (137) of Petersen and Pederson [14] to the second non-constant summand of (4.9). Both of those steps use the fact that $T$ and $\alpha$ are symmetric. The second equality of (4.10) follows from (511) of Petersen and Pederson [14], and the last from the definition of a Kronecker product.

The matrices $\frac{\partial \alpha}{\partial \alpha_{i,j}}$ for $1 \leq i \leq j \leq m$ in (4.10) do not depend on $\alpha$ since all the elements of each are zeros and ones. Furthermore, any matrix $\chi$ of dimension $m \times m$ can be expressed as a linear combination of them. By (4.10), the likelihood equations $\frac{\partial}{\partial \alpha_{i,j}} \mathcal{L}_1 = 0$ for $1 \leq i \leq j \leq m$ then imply that

$$tr\big((-(n-1)\alpha^{-1} + \alpha^{-1}YT^{-1}Y^T\alpha^{-1})\chi\big) = 0$$

for any matrix $\chi$ and, in particular, for $\chi = (-(n-1)\alpha^{-1} + \alpha^{-1}YT^{-1}Y^T\alpha^{-1})^T$. That is true if and only if $-(n-1)\alpha^{-1} + \alpha^{-1}YT^{-1}Y^T\alpha^{-1} = 0$, and we conclude that the likelihood equations are uniquely satisfied by $\alpha = \hat{\alpha}$ where

$$\hat{\alpha} = (n-1)^{-1}YT^{-1}Y^T, \tag{4.11}$$

or, equivalently, where (4.2) holds.

Applying (124), (125), and (137) of Peterson and Pedersen [14] to the expression on the right-hand side of (4.10) and substituting $\hat{\alpha}$ from (4.11) results in

$$\frac{\partial^2}{\partial \alpha_{i,j} \partial \alpha_{k,l}} \mathcal{L}_1(\hat{\alpha}) = -\frac{1}{2}(n-1)tr\left(\left(\hat{\alpha}^{-1}\frac{\partial \alpha}{\partial \alpha_{i,j}}\right)\left(\hat{\alpha}^{-1}\frac{\partial \alpha}{\partial \alpha_{k,l}}\right)\right) \tag{4.12}$$

for $1 \leq i \leq j \leq m$ and $1 \leq k \leq l \leq m$. The matrices $\hat{\alpha}^{-1}(\partial \alpha/\partial \alpha_{i,j})$ for $1 \leq i \leq j \leq m$ in (4.12) form a basis that spans a vector space, and the trace in (4.12) is as an inner product of two elements from that vector space. Consequently, the likelihood function's Hessian matrix evaluated at $\hat{\alpha}$ is the negative of a Grammian and therefore negative definite. This shows that $\mathcal{L}_1(\hat{\alpha})$ is indeed a maximum, so that $\hat{\alpha}$ in (4.2) is the unique maximum likelihood estimator for $p_1(x_1, x_2, \dots, x_{n-1}|x_n)$ in (4.4)

By (4.9), the log-likelihood function $\mathcal{L}_1$ depends on $Y$ and hence on the sample vectors $x_1, x_2, \dots, x_n$ only through $\hat{\alpha}$ in (4.11), or, equivalent, in (4.2). The sufficiency of $\hat{\alpha}$ then follows from the Fisher-Neyman theorem; see Halmos and Savage [18] for a formulation of that theorem covering the current case of dependent random variables. ∎

Returning to the assumptions of Theorem 3, let

$$p_1(x_1, x_2, \dots, x_{n-1}|x_n) \equiv \frac{P(X(t_1) \in dx_1, \dots, X(t_{n-1}) \in dx_{n-1}|X(t_n) = x_n)}{dx_1 \dots dx_{n-1}} \tag{4.13}$$

$$p_2(x_n) \equiv \frac{P(X(t_n) \in dx_n)}{dx_n} \tag{4.14}$$

and

$$p(x_1, x_2, \dots, x_n) \equiv \frac{P(X(t_1) \in dx_1, \dots, X(t_n) \in dx_n)}{dx_1 \dots dx_n} \tag{4.15}$$

They are related by

$$p(x_1, x_2, \dots, x_n) = p_1(x_1, x_2, \dots, x_{n-1}|x_n)p_2(x_n). \tag{4.16}$$

Following the usual convention, we will say that estimators are jointly sufficient if they are jointly sufficient with respect to (4.15) and are the maximum likelihood estimators if they maximize (4.15).

By Lemma 2, (4.13) coincides with the density of the same name in (4.4). Let $\Sigma \equiv \Gamma(t_n, t_n)$. If $Rank(\Sigma) = m$, then (4.14) has the familiar form,

$$p_2(x_n) = (2\pi)^{-m/2}|\Sigma|^{-1/2}exp\left(-\frac{1}{2}tr(\Sigma^{-1}x_n x_n^T)\right). \tag{4.17}$$

The density in (4.17) does not exist when $Rank(\Sigma) < m$, but, in that case, there exists a matrix $N$ with dimension $m \times (m - Rank(\Sigma))$ such that $Rank(N) = m - Rank(\Sigma)$ and $N^T \Sigma = 0$. Then, $N^T X(t_n)$ has a variance of zero and is therefore equal to its mean with probability one, i.e., $N^T X(t_n) = N^T E X(t_n) = 0$ with probability one. In this case, Section 8a.4 on pages 527-528 of Rao [16] shows that $X(t_n)$ will have the density,

$$p_2(x_n) = \begin{cases} (2\pi)^{-Rank(\Sigma)/2} |\Sigma|_+^{-1/2} exp\left(-\frac{1}{2} tr(\Sigma^- x_n x_n^T)\right), & N^T x_n = 0, \\ 0, & othewise. \end{cases} \quad (4.18)$$

where $|\Sigma|_+$ denotes the product of the non-zero eigenvalues of $\Sigma$ (and where $\Sigma^-$ denotes any generalized inverse of $\Sigma$ as before). When $Rank(\Sigma) = m$, the densities defined by (4.17) and (4.18) agree, so that (4.18) holds for the general case. We conclude using (4.16) that the density in (4.15) is equal to the product of (4.18) and $p_1(x_1, x_2, \ldots, x_{n-1} | x_n)$ from Lemma 3.

*Proof of Theorem 3:* When (4.3) holds, we can substitute $x_n x_n^T = t_n \hat{\alpha} - t_n^2 \hat{\beta}$ into (4.18) – noting that $N^T x_n = 0$ if and only if $0 = N^T x_n x_n^T N = N^T (t_n \hat{\alpha} - t_n^2 \hat{\beta}) N$ -- to show that (4.18) depends on the sample vector $x_n$ only through through $\hat{\alpha}$ and $\hat{\beta}$. The proof of Lemma 3 shows that $p_1(x_1, x_2, \ldots, x_{n-1} | x_n)$ depends on the sample vectors $x_1, x_2, \ldots, x_n$ only through $\hat{\alpha}$. Since both terms on the right-hand side of (4.16) depend on the sample vectors only through $\hat{\alpha}$ and $\hat{\beta}$, their joint sufficiency follows from the Fisher-Neyman theorem for any given $Rank(\Sigma)$.

Section 8a.5, page 532 of Rao [16] implies for the current case that

$$\hat{\Sigma} = x_n x_n^T \quad (4.19)$$

is the maximum likelihood estimator of $\Sigma \equiv \Gamma(t_n, t_n)$ in (4.18). In that case, (4.1) implies that any estimators $\hat{\alpha}$ and $\hat{\beta}$ satisfying $x_n x_n^T = t_n(\hat{\alpha} - \hat{\beta} t_n)$ maximize (4.18). In particular, (4.2) and (4.3) do. By Lemma 3, $\hat{\alpha}$ is the maximum likelihood estimator of (4.13), which does not depend on $\beta$. Since $\hat{\alpha}$ and $\hat{\beta}$ maximize the terms on the right-hand side of (4.16) individually, they are the maximum likelihood estimators as defined earlier.

By (4.1) and (4.2),

$$E\hat{\alpha} = \frac{1}{n-1}E\sum_{i=1}^{n-1}\frac{t_{i+1}^2 x_i x_i^T - t_i t_{i+1} x_i x_{i+1}^T - t_i t_{i+1} x_{i+1} x_i^T + t_i^2 x_{i+1} x_{i+1}^T}{t_i t_{i+1}(t_{i+1}-t_i)}$$

$$= \frac{1}{n-1}\sum_{i=1}^{n-1}\frac{t_{i+1}^2 \Gamma(t_i,t_i) - t_i t_{i+1}\Gamma(t_i,t_{i+1}) - t_i t_{i+1}\Gamma(t_{i+1},t_i) + t_i^2 \Gamma(t_{i+1},t_{i+1})}{t_i t_{i+1}(t_{i+1}-t_i)}$$

$$= \frac{1}{n-1}\sum_{i=1}^{n-1}\alpha = \alpha,$$

so that $\hat{\alpha}$ is unbiased. Then, by (4.1) and (4.3),

$$E\left(t_n(\hat{\alpha}-\hat{\beta}t_n)\right) = Ex_n x_n^T = \Gamma(t_n,t_n) = t_n(\alpha-\beta t_n)$$

so that $\hat{\beta}$ is unbiased as well. ∎

## 4.2   Posterior Moments

Our final result establishes a link between Gaussian process regression and linear regression when the covariance kernel defined by the assumptions of Theorem 3 is used for the former.

*Theorem 4: If $\left\{\hat{X}(t): 0 \leq t \leq T_0 < T \leq \infty\right\}$ is a zero-mean Gaussian process of dimension $m \times 1$ and (i) $\hat{X}(0) = 0$, (ii) samples $\hat{X}(t_i) = x_i$ for $0 < t_1 < t_2 < \ldots < t_n \leq T_0$ are known, and (iii) its covariance kernel $\hat{\Gamma}(\cdot,\cdot)$ is given by (4.1) with $\alpha = \hat{\alpha}$ from (4.2) and $\beta = \hat{\beta}$ from (4.3), then*

$$E\big[\hat{X}(t)\big|\hat{X}(t_1) = x_1, \ldots, \hat{X}(t_n) = x_n\big] = x_n + \left(\frac{x_n}{t_n} - (x_n^T x_n)^{-1}\hat{\alpha}x_n\right)(t-t_n) \tag{4.20}$$

and

$$E\big[\hat{X}(t_{n+p})\hat{X}^T(t_{n+p})\big|\hat{X}(t_1) = x_1, \ldots, \hat{X}(t_n) = x_n\big] = \hat{\Gamma}(t,t) - (x_n^T x_n)^{-2}\hat{\Gamma}(t,t_n)x_n x_n^T\hat{\Gamma}(t_n,t) \tag{4.21}$$

for $t_n \leq t \leq T_0$.

*Proof:* The kernel $\hat{\Gamma}(\cdot,\cdot)$ satisfies the conditions in Theorem 2 of Beutler [10] for the wide-sense Markov property. Since the Markov and wide-sense properties are equivalent for Gaussian processes, we can conclude that $\hat{X}$ is also a Markov process, even without the stronger assumptions of Corollaries 3 or 4. This implies that the conditional distribution of $\hat{X}(t)$ given $\hat{X}(t_i) = x_i$ for $t_1 < t_2 < \ldots < t_n \leq t$ is equal to the conditional distribution of $\hat{X}(t)$ given $\hat{X}(t_n) = x_n$. Then,

$$E\big[\hat{X}(t)\big|\hat{X}(t_1) = x_1, \ldots, \hat{X}(t_n) = x_n\big] = \hat{\Gamma}(t,t_n)\hat{\Gamma}(t_n,t_n)^- x_n \tag{4.22}$$

and

$$E\big[\hat{X}(t_{n+p})\hat{X}^T(t_{n+p})\big|\hat{X}(t_1) = x_1, \ldots, \hat{X}(t_n) = x_n\big] = \hat{\Gamma}(t,t) - \hat{\Gamma}(t,t_n)\hat{\Gamma}(t_n,t_n)^-\hat{\Gamma}(t_n,t) \tag{4.23}$$

for $t_n \leq t \leq T_0$ as follow from the formulas from the proof of Corollary 2. Since (4.22) and (4.23) hold for any generalized inverse $\hat{\Gamma}(t_n,t_n)^-$, they hold in particular for

$$\hat{\Gamma}(t_n,t_n)^- = \hat{\Gamma}(t_n,t_n)^+ = x_n(x_n^T x_n)^{-2} x_n^T, \tag{4.24}$$

where $a^+$ denotes the Moore-Penrose pseudo inverse; see Section 3.6 of Petersen and Pederson [14] for background. The last equality in (4.24) is easily verified from the definition of the Moore-Penrose pseudo inverse given that $\hat{\Gamma}_0(t_n,t_n) = x_n x_n^T$. The results in (4.20) and (4.21) immediately follow from (4.22)-(4.24). ∎